%% file: positivity.tex
\documentclass[11pt]{article}

\usepackage{epsfig}
\usepackage{array}
\usepackage{amsmath}
\usepackage{amssymb}
\usepackage{theorem}
\usepackage{url}

\makeatletter
\renewcommand{\subsection}[1]{\@startsection{subsection}{2}{0pt}%
{\medskipamount}{-4.5 pt}{\textbf}{#1}\hskip-\parindent\textbf{. }}
\makeatother

\newsavebox{\attributebox}

\theoremstyle{change}
\theorempreskipamount=\bigskipamount
\theorempostskipamount=\smallskipamount

\newtheorem{cor}[subsection]{Corollary. --- }

\newtheorem{prop}[subsection]{Proposition. --- }
\newtheorem{thm}[subsection]{Theorem. --- }

\theorembodyfont{\rmfamily}

\newtheorem{ex}[subsection]{Example. --- }

\input NormalForm_ntreemacros

\newcommand{\proof}{\noindent\textit{Proof. --- }}
\newcommand{\PBS}[1]{\let\temp=\\#1\let\\=\temp}

\newcommand{\1}{^{-1}}
\newcommand{\<}{\mathopen{<}}
\renewcommand{\>}{\mathclose{>}}
\renewcommand{\b}{v+v^{-1}}
\newcommand{\cb}[3]{c_{\<#3,#1,#2]}}
\newcommand{\cf}[3]{c_{[#1,#2,#3\>}}
\renewcommand{\H}{{\cal H}}
\newcommand{\Hfour}{H${}_4$}
\newcommand{\Hthree}{H${}_3$}

\newcommand{\LR}{{\rm LR}}

\newcommand{\Z}{{\bf Z}}

\begin{document}

\title{Positivity results for the Hecke algebras
of non-crystallographic finite Coxeter groups}
\author{Fokko du Cloux}

\maketitle

\begin{flushleft}
Institut Camille Jordan\footnote{The author gratefully acknowledges the
hospitality and technical assistance of the Centre de Calcul Medicis,
Laboratoire STIX, FRE 2341 CNRS, and of Maply, UMR 5528 CNRS, where the
computations for this paper have been carried out.}\\
UMR 5208 CNRS\\
Universit\'e Lyon-I\\
69622 Villeurbanne Cedex FRANCE\\
{\tt ducloux@math.univ-lyon1.fr}
\end{flushleft}

\noindent{\em Abstract.} This paper is a report on a computer check of
some important positivity properties of the Hecke algebra in type \Hfour, 
including the non-negativity of the structure constants in the Kazhdan-Lusztig 
basis. This answers a long-standing question of Lusztig's. The same
algorithm, carried out by hand, also allows us to deal with the case of 
dihedral Coxeter groups.

\section{Statement of the problem}\label{section:introduction}

\subsection{}
Let $W$ be a Coxeter group, $S$ its set of distinguished generators, and denote
$\leq$ the Bruhat ordering on $W$. Denote $\H$ the Hecke algebra of W over
the ring of Laurent polynomials $A=\Z[v,v\1]$, where $v$ is an indeterminate.
We refer to \cite{geck_pfeiffer:2000} for basic results about Coxeter groups
and Hecke algebras; we just recall here that $\H$ is a free $A$-module with
basis $(t_y)_{y\in W}$, where the algebra structure is the unique one which
satisfies
$$
t_s.t_y=\begin{cases}t_{sy}&\text{if $sy>y$}\\
(v-v\1)t_y+t_{sy}&\text{if $sy<y$}\end{cases}
$$
(here we use $t_y=v^{-l(y)}T_y$, where the $T_y$ satisfy the more familiar
relation $T_s.T_y=(q-1)T_y+qT_{sy}$ when $sy<y$, with $q=v^2$.) Then there is
a unique ring involution $i$ on $\H$ such that $i(v)=v\1$, and 
$i(t_y)=t_{y\1}\1$.

The Kazhdan-Lusztig basis of $\H$ is the unique family $(c_y)_{y\in W}$
such that {\em (a)} $i(c_y)=c_y$ and {\em (b)} 
$c_y=t_y+\sum_{x<y}p_{x,y}t_x$, with $p_{x,y}\in v\1\Z[v\1]$; in particular
we find that $c_s=t_s+v\1$ for all $s\in S$. It turns out
that $P_{x,y}=v^{l(y)-l(x)}p_{x,y}$ is a polynomial in $q$, the {\em
Kazhdan-Lusztig polynomial} for the pair $x,y$. For any pair $(x,y)$ of
elements in $W$, write
$$
c_x.c_y=\sum_{z\in W}h_{x,y,z}c_z
$$
(in other words, the $h_{x,y,z}\in A$ are the structure constants of the Hecke
algebra in the Kazhdan-Lusztig basis).

A number of the deeper results in the theory of Hecke algebras depend on
positivity properties of the polynomials $P_{x,y}$ and $h_{x,y,z}$. More
precisely, consider the following properties~:

\begin{itemize}
\item[$P_1$:]the polynomials $P_{x,y}$ have non-negative coefficents;
\item[$P_2$:]the $P_{x,y}$ are decreasing for fixed $y$, in the sense that
if $x\leq z\leq y$ in $W$, the polynomial $P_{x,y}-P_{z,y}$ has non-negative
coefficents;
\item[$P_3$:]the polynomials $h_{x,y,z}$ have non-negative coefficients.
\end{itemize}

Properties $P_1$ and $P_3$ are basic tools in the study of Kazhdan-Lusztig
cells and the asymptotic Hecke algebra; they have been proved in
\cite{kl:1980} and \cite{springer:1982} for crystallographic $W$ using deep 
properties of intersection cohomology (see the remarks in section 3 of Lusztig 
\cite{lusztig:1985a} for the case where $W$ is infinite.) Property $P_2$ is 
proved for finite Weyl groups in \cite{irving:1988}, using the description of 
Kazhdan-Lusztig polynomials in terms of filtrations of Verma modules.

None of these geometric or representation-theoretic interpretations are 
available in the non-crystallographic case, and the validity in general 
of the positivity properties above are among the main open problems in the 
theory of Coxeter groups. Let us concentrate on the case where $W$ is finite. 
It is easy to see that for the validity of $P_1$--$P_3$ we may reduce to the
case where $W$ is irreducible.

\subsection{}
The case where $W$ is dihedral is simple enough to be carried out by hand;
we have included the computation of the $h_{x,y,z}$ in section
\ref{section:dihedral}. This leaves us only with the two groups \Hthree\ and 
\Hfour, and since the former is contained in the latter, the only case we need 
to consider is the Coxeter group of type \Hfour, of order $14400$, with Coxeter
diagram
$$
\CoxeterDiagramHfour
$$
Despite its rather modest size, the group \Hfour\ poses a redoubtable
computational challenge, even for present-day computers. It is all the more
remarkable that property $P_1$ was checked already in 1987 by Dean Alvis 
\cite{alvis:1987} by explicitly computing all the Kazhdan-Lusztig polynomials.
Quite a feat with the hardware of the time! Unfortunately, Alvis's programs 
have never been made available; to my knowledge, the only publicly available 
computer program capable of carrying out this computation is my own program 
{\tt Coxeter} \cite{du_cloux:coxeter}, which does it in about one minute on
a modern-day personal computer. This still leaves open properties $P_2$ and 
$P_3$; the main purpose of this paper is to report on a computation, carried 
out as one of the first applications of version 3 of {\tt Coxeter}, by which 
we prove~:

\begin{thm}\label{thm:positivity}
Properties $P_2$ and $P_3$ hold for the Coxeter group of type \Hfour.
\end{thm}

\begin{cor}\label{cor:positivity}
The fifteen conjectures labelled \textbf{P1--P15} in Chapter 14 of Lusztig 
\cite{lusztig:2003}
hold for the Hecke algebra of the Coxeter group of type \Hfour, and in fact
for the equal parameter Hecke algebra of any finite Coxeter group.
\end{cor}

\proof Lusztig shows in Chapter 15 of \cite{lusztig:2003} that in the equal
parameter case (which is automatic in type \Hfour), the fifteen conjectures all
follow from $P_1$ and $P_3$. In view of our earlier remarks
on the reduction to the irreducible case, and of the dihedral computation in
section \ref{section:dihedral}, the same argument can be applied to any finite
Coxeter group.

\section{Methodology}\label{section:methodology}

\subsection{}The verification of $P_1$ and $P_2$ is straightforward~: one
simply runs through the computation of the $P_{x,y}$ for all $x\leq y$ in
the group. In fact, from the well-known property $P_{x,y}=P_{sx,y}$ whenever
$sx>x$, $sy<y$ (cf.\ \cite{kl:1979} (2.3.g)), and the analogous property on 
the right, for $P_1$ it suffices to
consider the cases where $\LR(x)\supset \LR(y)$, where we denote 
$L(x)=\{s\in S\mid sx<x\}$ (resp.\ $R(x)=\{s\in S\mid xs<x\}$) the left
(resp.\ right) descent set of $x$, and $LR(x)=L(x)\coprod R(x)\subset S\coprod 
S$; we call such pairs $(x,y)$ {\em extremal pairs.} Taking also into
account the symmetry $P_{x,y}=P_{x\1,y\1}$, there are $2\,348\,942$ cases to 
consider, which are easily tabulated by the program. The non-negativity of the 
polynomials is checked as they are found. For $P_2$, one also easily reduces
to extremal pairs $(x,y)$ and $(z,y)$. The tough computation is for $P_3$; here
there are {\it a priori} $14\,400^3=2\,985\,984\,000\,000$ (almost three
trillion!) polynomials $h_{x,y,z}$ to be computed, and the only obvious
symmetry is $h_{x,y,z}=h_{y\1,x\1,z\1}$ (but, as explained below, we don't
even use that.)

\subsection{}\label{subsection:algorithm}
The algorithm used in the computation is simple. For a fixed $y$, 
we compute the various $c_x.c_y$ by induction on the length of 
$x$, starting with $c_e.c_y=c_y$, where $e$ denotes the identity element of 
$W$. To carry out the induction, we choose
any generator $s\in S$ such that $sx<x$, and write~:
$$
c_x=c_s.c_{sx}-\underset{sz<z}{\sum_{z<sx}}\;\mu(z,sx)c_z\eqno(1)
$$
where as usual $\mu(x,y)$ denotes the coefficient of $v\1$ in $p_{x,y}$ (which
is also the coefficient of degree $\frac{1}{2}(l(y)-l(x)-1)$ in $P_{x,y}$, and
in particular is zero when the length difference of $x$ and $y$ is even.)

Then we may assume that $c_{sx}.c_y$ is already known, and similarly for the
various $c_z.c_y$, so we are reduced to multiplications of the form $c_s.c_u$,
for $s\in S$ and $u\in W$. When $su>u$ this is read off from formula $(1)$ 
with $x=su$; and when $su<u$ one simply has $c_s.c_u=(v+v\1)c_u$ (see for
instance \cite{kl:1979}).

The information that is required for this computation is encoded in the
$W$-graph of the group; once this graph is known, everything else just involves
elementary operations on polynomials. Actually, it is obvious that the
$h_{x,y,z}$ are in fact polynomials in $(v+v\1)$, so they are determined by
their positive-degree part; it is this part that we compute and keep in memory.
The only complication arises from the fact that we need to be careful about
integer overflow; it turns out in fact that for \Hfour\ all the coefficients
of the $h_{x,y,z}$ fit into a $32$-bit unsigned (and even signed) integer, but 
not by all that much: the largest coefficient that occurs is $710\;904\;968$, 
which is only about six times smaller than $2^{32}=4\;294\;967\;296$.

\subsection{}From the computational standpoint, this procedure has a number
of desirable features. First and foremost, once the $W$-graph of the group has
been determined, the problem splits up into $14400$ independent
computations, one for each $y$, so we can forget about the computation for
a given $y$ when passing to the next (this would not be true if we tried to
use the symmetry $h_{x,y,z}=h_{y\1,x\1,z\1}$.) This is advantageous in terms
of memory usage, could be used to parallelize the computation if necessary
(it turns out that it hasn't been), and also means that the computation can
be harmlessly interrupted (either voluntarily or involuntarily), at least if
its progress is recorded somewhere~: basically, the only penalty to pay for 
picking up an interrupted computation is the recomputation of the $W$-graph,
which takes about half a minute. This is very valuable for computations
running over several days, where there is always the risk of a system crash
or power failure.

On the other hand, it is essential that for a fixed $y$, the full table of
$h_{x,y,z}$ be stored in memory. In practice, there are many repetitions
among these polynomials; so we store them in the form of a table of
$14400^2=207\,360\,000$ pointers. Initially, this is the main memory 
requirement
of the program; it is interesting to note that the cost doubles, from about
$800$ MB to about $1.6$ GB, when we pass from a 32-bit to a 64-bit computer.
It turns out that for a fixed $y$, the additional memory required to store
the actual polynomials is small, and never exceeds $300$ MB. So the full 
computation runs confortably in $2$ GB of memory, and barely exceeds $1$ GB 
on a 32-bit machine.

It is much more difficult to try to write down a full table of all the
polynomials that occur as $h_{x,y,z}$. I have done this a number of times,
but with the memory available on the machines to which I have had access,
it has been necessary to split up the computation in about a hundred pieces
to avoid memory overflow in the polynomial store, to keep the corresponding
files in compressed form to avoid overflowing the hard disk, and then to
merge those one hundred compressed files into a single compressed file,
eliminating repetitions.
At some point I have needed to store about $30$ GB of compressed files---not 
something administrators are very happy about! In view of these difficulties,
I have chosen not to make that version of the program available for the time
being.

\subsection{}
The computation has been done a number of times (writing files of all the
distinct polynomials): at the Ecole Polytechnique,
Centre de Calcul M\'edicis, Laboratoire STIX, FRE 2341 CNRS,
on several computers, including a Compaq Alpha EV68 and an AMD Opteron server
with $4$ and $8$ GB of main memory respectively, where it has taken about 5 
days of CPU time, and twice at the Universit\'e Lyon I, Maply, UMR 5585 CNRS,
on a Xeon processor with $4$ GB of main memory, where it took 80 hours (not
counting the final file merging pass, which took another 80 hours or so.) The
program as presented here has run on our 2.7 GHz AMD Opteron server at the 
Institut Camille Jordan, using less than 2 GB of memory, in a little under 85 
hours.

On the technological side, it seems that the time was just right for this
computation~: it makes full use of the 3Ghz processors, at least 2GB central
memories, and 100+ GB hard disks that are found on typical low-end servers
today. It would still be beyond the grasp of most present-day personal 
computers, however, although that, too, is changing fast!

\section{Verification}

\subsection{}Let's come now to the thorny issue of verification~: what is the
amount of trust that can be put in a result like this~? An obvious prerequisite
is the availablility of the source code of the program that carries out the
computation; this may be downloaded at
\url{http://igd.univ-lyon1.fr/~ducloux/coxeter/coxeter3/positivity}.

This is actually the source code of an especially modified version of
{\tt Coxeter3}. All the extra code is contained in the file {\tt special.cpp};
all the other files are identical to the ones in {\tt Coxeter3}.

\subsection{}\label{subsection:commands}
In addition to those already available in {\tt Coxeter3},
the following commands are defined~:

\begin{itemize}\itemsep-4pt
\item[--]{\tt klplist}~: prints out a list of all the distinct Kazhdan-Lusztig
polynomials which occur in the group (so that in particular, one may re-check
property $P_1$, although this was done already by the computer in the course
of the computation);
\item[--]{\tt decrklpol}~: checks property $P_2$;
\item[--]{\tt positivity}~: checks the non-negativity of the $h_{x,y,z}$;
\item[--]{\tt cycltable}~: prints out a table of the $c_x.c_y$ for a fixed $y$;
\item[--]{\tt cprod}~: prints out a single product $c_x.c_y$;
\end{itemize}

For type $H_4$, on a decent server, the {\tt cycltable} command should not 
take more than two minutes; {\tt cprod} should usually take less than a minute
(Note that the first call will take longer than subsequent ones, because the
$W$-graph must be computed the first time.) So these commands 
give local access to the multiplication table of the Hecke algebra, thus
opening up the ``black box'' a little.

The {\tt positivity} command for type $H_4$ can also be executed through the
little stand-alone command {\tt coxbatch} that I have included; this will run
the computation in the background. It writes any errors in {\tt error\_log}, 
and records the progress of the computation in {\tt positivity\_log}. After a
successful run, {\tt error\_log} should be empty, and {\tt positivity\_log}
should end with the line~:
$$
\texttt{14399: maxcoeff = 710904968}
$$
(elements of the group are numbered from $0$ to $14399$.)

\subsection{}As was explained in \ref{subsection:algorithm}, the computation is
entirely elementary once the $W$-graph of the group is known.
The trust that one may place in this ingredient should be rather high, in my
opinion, as it is computed with an algorithm which in simpler cases has been
checked against other programs, and which even for \Hfour\ has been checked
against other algorithms for the same computation. (The latter check is perhaps
more convincing than the former, for it may very well happen that some of the 
nastiest configurations occur for the first time in type \Hfour, or even
only in type \Hfour, as it is such an exceptional group.)

For the actual non-negativity check, we are as far as I know in entirely 
uncharted territory. However, the code for the computation of the $h_{x,y,z}$
from the $W$-graph is really quite simple, so it can 
rather convincingly be checked by inspection.
Another check is as follows~: the order in which the computations are performed
depends on the choice of a descent generator for $x$. By default, we always
choose the first such generator in the internal numbering of {\tt Coxeter};
however, it is easy enough to replace this by other ``descent strategies''
(for instance, choosing the first generator in some other ordering.) This
will lead to a very different flow of recursion. I have done this replacing
the first descent generator by the last one, and obtained the same output
file, which is as it should be.

\section{The case of dihedral groups}\label{section:dihedral}

\subsection{} Let us now show how the algorithm described in 
\ref{subsection:algorithm} can be carried out by hand in the case where $W$ is 
dihedral. This is not difficult and may well be known to experts, but I 
haven't been able to find the results in the literature.
In Chapter 17 of \cite{lusztig:2003}, Lusztig computes the $h_{x,y,z}$ for
an infinite dihedral group with {\em unequal} parameters. I have
not been able to determine if his formulas can be specialized to yield
the statement in Proposition \ref{prop:infinite dihedral} (although of course
our computations for the infinte case are very similar to (and much simpler 
than) those in \cite{lusztig:2003}.)

\subsection{}\label{subsection:dihedral notation}
Assume first that $W$ is infinite dihedral. Let 
$S=\{s_1,s_2\}$, and for each $i\geq 0$ denote
\begin{align*}
&[1,2,i\>=s_1s_2\ldots\qquad[2,1,i\>=s_2s_1\ldots\\
&\<i,1,2]=\ldots s_1s_2\qquad\<i,2,1]=\ldots s_2s_1
\end{align*}
where in each case there are $i$ terms in the product. When we need 
indeterminate generators, we will also use notation such as $[s,t,i\>$ for
$s\neq t$ in $S$. Denote for simplicity
$c_1=c_{s_1}$, $c_2=c_{s_2}$. It is well-known that for dihedral groups all
the Kazhdan-Lusztig polynomials $P_{x,y}$ are equal to $1$ for $x\leq y$;
incidentally, this proves that $P_1$ and $P_2$ are trivially verified. It
follows that formula \ref{subsection:algorithm} (1) reduces to
\begin{align*}
&c_1.c_2=c_{s_1.s_2}=\cf122\\
&c_1.\cf 21i = \cf12{i+1}+\cf12{i-1}\quad\text{for $i>1$}
\end{align*}
and of course $c_1.\cf 12i = (\b)\;\cf12i$ for $i>0$, and similar formul\ae\
for the left multiplication by $c_2$.

\subsection{}\label{subsection:recursion}
Now fix an element of the Kazhdan-Lusztig basis, that without
loss of generality we may assume to be of the form $\cb12k$, with $k>0$. 
Let $s$ be the first generator in $\<1,2,k]$, so that $\cb 12k=\cf stk$, and 
let $t\neq s$ be the other
element of $S$. We wish to compute the various $\cb sti.\cb12k$, 
$\cb tsi.\cb12k$, as $i$ varies. Note that since the $\cb 12j$, $j>0$, form a 
basis of a left cell representation of $\H$, it is {\it a priori} clear that 
these products will be $A$-linear combinations of the $\cb 12j$.
Start with the $\cb sti.\cb12k$, where we may assume $i>0$. For the first two 
values of $i$ we get
\begin{align*}
&\cb st1.\cb12k=c_t.\cf stk=\cf ts{k-1}+\cf ts{k+1}=\cb12{k-1}+\cb12{k+1}\\
&\cb st2.\cb12k=c_s.c_t.\cb12k=\cb12{k-2}+2\cb12k+\cb12{k+2}
\end{align*}
when $k>2$, and
\begin{align*}
&\left\{\begin{aligned}
&c_t.\cb 121=\cb 122\\
&c_s.c_t.c_2=\cb 121+\cb123
\end{aligned}
\qquad\text{for $k=1$}\right.\\
&\left\{\begin{aligned}
&c_t.\cb122=\cb211+\cb213\\
&c_s.c_t.\cb122=2\cb212+\cb214
\end{aligned}
\qquad\text{for $k=2$}\right.
\end{align*}
Now applying the procedure from \ref{subsection:algorithm} yields the recursion
formula:
$$
\cb sti.\cb12k=c_r.\cb st{i-1}.\cb12k-\cb st{i-2}.\cb12k\quad\text{for $i>2$}
$$
where $r\in\{s,t\}$ is the first term in $\cb sti$. It follows easily that
the non-zero terms in $\cb sti.\cb12k$ all correspond to indices $j$ of the
same parity, which changes when we go from one $i$ to the next.
Consequently, all these terms have a first generator {\it not} equal to $r$.
If we write the coefficients in rows, we see that the coefficient at position 
$j$ in row $i>2$ is obtained by adding the coefficients at positions $j-1$ and 
$j+1$ in row $i-1$, and subtracting the coefficient at position $j$ in row 
$i-2$ . For example, when $k=3$, the table looks like this (with dots 
indicating zeroes)~:
$$
\begin{matrix}
&k-2&k-1&k&k+1&k+2&k+3&k+4&k+5\\
\\
\text{$i=1$\qquad}&\cdot&1&\cdot&1&\cdot&\cdot&\cdot&\cdot\\
\text{$i=2$\qquad}&1&\cdot&2&\cdot&1&\cdot&\cdot&\cdot\\
\text{$i=3$\qquad}&\cdot&2&\cdot&2&\cdot&1&\cdot&\cdot\\
\text{$i=4$\qquad}&1&\cdot&2&\cdot&2&\cdot&1\\
\text{$i=5$\qquad}&\cdot&1&\cdot&2&\cdot&2&\cdot&1\\
\end{matrix}
$$
Due to the exceptional multiplication
rules $c_1.c_2=\cb212$ and $c_2.c_1=\cb122$, the first $1$ that would appear
in column 0 should be omitted (this occurs for $i = k$); another way of 
stating this is that we should
run the algorithm in the half-plane $j>0$, and omit undefined terms. In this
form, the procedure is valid for all values of $k>0$. Note that for $i>k$ the
rows have constant length $k+1$, the next one being just the previous one
shifted one unit to the right.

The case of the $\cb tsi.\cb12k$ is similar but simpler; it appears as
Proposition 7.7 {\it (a)} in Chapter 7 of \cite{lusztig:2003}. Here all the 
non-zero coefficients are equal to $\b$, and we have a second order recursion 
as previously. Again the coefficients start expanding pyramidally, there is a
``reflection'' at column $0$, and then we get a band of constant width $k$.

The final result may be stated as follows:

\begin{prop}\label{prop:infinite dihedral}
Let $W$ be an infinite dihedral group, and adopt the notation from
\ref{subsection:dihedral notation} Let $k>0$, and let $s$ be the first
generator in $\cb12k$ ({\it i.e.}, $s=1$ if $k$ is even, $s=2$ if $k$ is odd.)
Let $t$ be the other element in $S$. Then we have~:
\vspace{-\bigskipamount}
\begin{multline*}
\cb sti.\cb12k=\\\left\{
\begin{alignedat}{2}
&\cb12{k-i}+2\cb12{k-i+2}+\ldots+2\cb12{k+i-2}+\cb12{k+i}&&\quad(0<i<k)\\
&2\cb12{2}+\ldots+2\cb12{2k-2}+\cb12{2k}&&\quad(i=k)\\
&\cb12{i-k}+2\cb12{i-k+2}+\ldots+2\cb12{i+k-2}+\cb12{i+k}&&\quad(i>k)\\
\end{alignedat}\right.
\end{multline*}
(where for $k=1$ the entry for $i=k$ should be interpreted as $\cb122$), and 
similarly
\vspace{-\bigskipamount}
\begin{multline*}
\cb tsi.\cb12k=\\\left\{
\begin{alignedat}{2}
&(\b)\;\cb12{k-i+1}+\ldots+(\b)\;\cb12{k+i-1}&&\quad(0<i<=k)\\
&(\b)\;\cb12{i-k+1}+\ldots+(\b)\;\cb12{i+k-1}&&\quad(i>=k)\\
\end{alignedat}\right.
\end{multline*}
with of course similar formul\ae\ when $\cb12k$ is replaced by $\cb21k$
\end{prop}

\subsection{}\label{subsection:finite dihedral}
Now consider the case where $W$ is a finite dihedral group of
order $2m$, $m\geq2$. Then of course $\cb12i$ is defined only for $i\leq m$,
and moreover the action of $c_1$ and $c_2$ on $\cb12m=\cb21m$ is given by
$$
c_1.\cb12m=c_2.\cb12m=(\b)\cb12m
$$
The consequence is that the recursive pattern described above gets modified
starting from $i=m+1-k$. For that row, the expression
$\cb12{m-1}+\cb12{m+1}$ that would have been obtained by applying the
appropriate $c_s$ to $\cb12m$ should be replaced by $(\b)\cb12m$. The net 
effect is that the algorithm splits in two independent parts: one one hand
we run the same algorithm as for the infinite case, but this time within the
strip $0<j<m$; on the other hand, starting from $i=m-k$, we add a term of
the form $a_i\cb12m$, with $a_{m-k}=1$, $a_{m-k+1}=\b$, 
$a_{m-k+2}=v^2+2+v^{-2}$,
and 
$$
a_i=(\b)a_{1-1}-a_{i-2}=v^d+2v^{d-2}+\ldots+2v^{-d+2}+v^{-d}
$$
for $d=k+i-m>2$. The procedure goes on until $i=m$, at which point the row 
from the first algorithm
has disappeared altogether, and only the multiple of $\cb12m$ is left.
For example, when $m=9$, $k=6$, and the $\cb12i.\cb12k$ (we have $s=s_1$ in
this case), we get the following table for the first part of the algorithm~:
$$
\begin{matrix}
&k-5&k-4&k-3&k-2&k-1&k&k+1&k+2\\
\\
\text{$i=1$\qquad}&\cdot&\cdot&\cdot&\cdot&1&\cdot&1&\cdot\\
\text{$i=2$\qquad}&\cdot&\cdot&\cdot&1&\cdot&2&\cdot&1\\
\text{$i=3$\qquad}&\cdot&\cdot&1&\cdot&2&\cdot&2&\cdot\\
\text{$i=4$\qquad}&\cdot&1&\cdot&2&\cdot&2&\cdot&1\\
\text{$i=5$\qquad}&1&\cdot&2&\cdot&2&\cdot&1&\cdot\\
\text{$i=6$\qquad}&\cdot&2&\cdot&2&\cdot&1&\cdot&\cdot\\
\text{$i=7$\qquad}&1&\cdot&2&\cdot&1&\cdot&\cdot&\cdot\\
\text{$i=8$\qquad}&\cdot&1&\cdot&1&\cdot&\cdot&\cdot&\cdot\\
\end{matrix}
$$
and row $9$ is zero. The obvious symmetry in the shape is a general fact.

\medskip

Note that if we look at the $\Z$-basis of $\H$ afforded
by the $v^d c_w$, $d\in\Z$, and give $v^d\cb12j$ degree $j+d$, then the
sum of the coefficients of $\cb sti.\cb12k$ in each degree is the same for 
the finite and the
infinite cases; whatever is missing from the infinite picture is exactly
reflected in the coefficient of $\cb12m$. So a concise statement of the
result is as follows~:

\begin{prop}\label{prop:finite dihedral}
Let $W$ be dihedral of order $2m$, $m\geq2$. Then the formul\ae\
from Proposition \ref{prop:infinite dihedral} remain valid, except that
one must have $i\leq m$, and that any expression of the form 
$$
\cb12{m+d}+\cb12{m-d}\qquad d>0
$$
that can be taken out of the formula should be replaced by 
$(v^d+v^{-d})\cb12m$.
\end{prop}

\begin{ex}
Pursuing the earlier example where $m=9$ and $k=6$, and taking for instance
the case where $i=6$, putting together the sixth row in the above table
and the expression for the coefficient $a_i$, which corresponds to $d=3$,
we get:
\begin{multline}
\cb126.\cb126=\\
2\cb122+2\cb124+\cb126+(v^3+2v+2v^{-1}+v^{-3})\;\cb129
\end{multline}
The corresponding expression for the infinite group would be
$$
2\cb122+2\cb124+2\cb126+2\cb128+2\cb12{10}+\cb12{12}
$$
from which we recover $(1)$ by taking out out one copy of $\cb126+\cb12{12}$ 
and two copies of $\cb128+\cb12{10}$.

\medskip

It is easy to get many other examples from the program --- of course playing
with the program is how the above statements were found in the first place.
\end{ex}

\section{Questions}

\subsection{}\label{subsection:unimodality}On looking at the lists of 
polynomials which occur as $h_{x,y,z}$, one immediately notices that they
are not only non-negative, but have a much stronger positivity property~:
if we denote $d$ the degree of $h_{x,y,z}$, then $v^dh_{x,y,z}$ is a
polynomial in $q=v^2$, which is {\em unimodal} (recall that this means
that the coefficients increase up to a point, which in this case has to be
the middle because of the symmetry of the $h_{x,y,z}$, and decrease from 
there.) In the course of the computation, the program checks unimodality for 
all $h_{x,y,z}$, and prints an error message on \texttt{error\_log} in case of
failure. Hence unimodality holds for the Hecke algebra of type \Hfour.

\subsection{}For Weyl groups, there is one case where it is easy to prove that 
unimodality holds~: {\it viz.}
the case where $y=w_0$ is the longest element in the group. From the properties
of the $\leq_\LR$ preorder it is clear that $A.c_{w_0}$ is a two-sided ideal
in $\H$; so we may write $c_x.c_{w_0}=h_xc_{w_0}$, where $h_x=h_{x,w_0,w_0}$.
Now it is clear that $t_s.c_{w_0}=v.c_{w_0}$ for all $s$, hence 
$t_x.c_{w_0}=v^{l(x)}c_{w_0}$, and 
$$
c_x.c_{w_0}=\sum_{z\leq x}p_{z,x}v^{l(z)}c_{w_0}
$$
from which it follows immediately, using the expression of the intersection
homology Poincar\'e polynomial in terms of Kazhdan-Lusztig polynomials
(\cite{kl:1980} Theorem 4.3.) that $v^{l(x)}h_x$ is equal to the Poincar\'e
polynomial of the intersection (hyper)cohomology of the Schubert variety
$X_x$. The unimodality then follows from the so-called hard Lefschetz theorem.
As far as I know, the unimodality property for general $h_{x,y,z}$ is an
open question, even for Weyl groups.

\subsection{}Clearly, all the results about the Hecke algebra of type \Hfour\
which are stated in this paper point to the fact that there is a hidden 
geometry here that is begging to be discovered. Hopefully, the facts
about this geometry which the program opens up will help us understand what
is going on, and serve as a guide towards the solution. I should be very
happy if this turns out to be the case.

\bigskip

\noindent{\bf Acknowledgement}

\medskip

I am grateful to the referee for some helpful remarks, and for pointing out
that the dihedral case is {\em not} trivial!

\bibliographystyle{plain}
\bibliography{positivity}

\end{document}

%% file: NormalForm_ntreemacros
\newcount\edgeno
\newdimen\edgelength
\newbox\chainbox
\newskip\leftchainskip

\def\ulap#1{\vbox to 0pt{\vss #1}}
\def\dlap#1{\vbox to 0pt{#1\vss}}
\def\udlap#1{\vbox to 0pt{\vss #1\vss}}
\def\rllap#1{\hbox to 0pt{\hss #1\hss}}

\def\smallbox#1{\udlap{\rllap{#1}}}

\def\smalldots{\smallbox{$\ldots$}}
\def\smallcirc{\smallbox{$\circ$}}
\def\smallnum#1{\smallbox{$\scriptscriptstyle #1$}}

\def\hgap{\hskip.1 cm}
\def\bighgap{\hskip.15 cm}
\def\vgap{\vskip.1 cm}
\def\bigvgap{\vskip.15 cm}

\def\vgapneg{\hskip -.1 cm}

\def\putbox#1#2#3{\vbox{\vskip#2\edgelength\hbox{%
\hskip#1\edgelength #3\hskip-#1\edgelength}%
\vskip-#2\edgelength}}

\def\vertexlabel#1{\hbox{$\scriptscriptstyle #1$}}

\def\rhedge#1{\udlap{\hbox to\edgelength{%
\hskip .5\wd\chainbox\hgap\leaders\hrule\hfil\hgap%
\global\setbox\chainbox\vertexlabel{#1}%
\hskip.5\wd\chainbox\smallnum{#1}}}}

\def\lhedge#1{\udlap{\hbox to \edgelength{%
\leftchainskip.5\wd\chainbox%
\global\setbox\chainbox\vertexlabel{#1}%
\smallnum{#1}\hskip .5\wd\chainbox%
\hgap\leaders\hrule\hfilneg\hgap%
\hskip\leftchainskip}}}

\def\dvedge#1{%
\rllap{\vbox to\edgelength{\offinterlineskip%
\vskip .5\ht\chainbox\vgap\leaders\vrule\vfil\vgap%
\global\setbox\chainbox\vertexlabel{#1}%
\vskip.5\ht\chainbox\smallnum{#1}}}}

\def\uvedge#1{%
\rllap{\vbox to\edgelength{\offinterlineskip%
\vskip -.5\ht\chainbox%
\vgapneg\leaders\vrule\vfil\vgapneg%
\global\setbox\chainbox\vertexlabel{#1}%
\vskip -.5\ht\chainbox\smallnum{#1}}}}

\def\shedge{\vbox to 0 pt{\vss\hbox to\edgelength{%
\bighgap\leaders\hrule\hfil\bighgap}\vss}}
\def\svedge{\hbox to 0 pt{\hss\vbox to\edgelength{%
\bigvgap\leaders\vrule\vfil\bigvgap}\hss}}

\def\dhedge{\udlap{\hbox to\edgelength{%
\hss\smalldots\hss}}}

\def\uclabel#1{\hbox to 0 pt{\hss\vbox to 0 pt%
{\vss\hbox{$\scriptstyle{#1}$}\vskip6 pt}\hss}}
\def\lclabel#1{\hbox to 0 pt{\hss\vbox to 0 pt{%
\vss\hbox{$\scriptstyle{#1}$\hskip6 pt}\vss}}}
\def\uelabel#1{\hskip-\edgelength\hbox to\edgelength%
{\hss\vbox to 0 pt%
{\vss\hbox{$\scriptstyle{#1}$}\vskip4 pt}\hss}}
\def\delabel#1{\hskip-\edgelength\hbox to\edgelength%
{\hss\vbox to 0 pt%
{\vskip4 pt\hbox{$\scriptstyle{#1}$}\vss}\hss}}

\def\dlabel#1{\hbox to -\edgelength{\hss\dlap%
{\vskip4 pt\vertexlabel{#1}}\hss}}
\def\llabel#1{\vbox to -\edgelength
{\vss\llap{\vertexlabel{#1}\hskip4 pt}\vss}}

\def\ulabel#1{\hbox to -\edgelength{\hss\ulap%
{\vertexlabel{#1}\vskip4 pt}\hss}}

\def\posulabel#1{\hbox to \edgelength{\hss\ulap%
{\smallnum{#1}\vskip4 pt}\hss}}

\def\dtedge#1#2{\dlap{\offinterlineskip%
\dvedge{#2}\llabel{#1}}}
\def\ltedge#1#2{\llap{\ulabel{#1}\lhedge{#2}}}
\def\rtedge#1#2{\rlap{\rhedge{#2}\ulabel{#1}}}

\def\de#1#2#3#4{\putbox{#1}{#2}{\dtedge{#3}{#4}}}

\def\le#1#2#3#4{\putbox{#1}{#2}{\ltedge{#3}{#4}}}
\def\re#1#2#3#4{\putbox{#1}{#2}{\rtedge{#3}{#4}}}

\def\hchain#1{\vbox{\edgeno0%
\halign{&\re{\the\edgeno}0##\global\advance\edgeno 
by 1\cr #1\cr}}}
\def\vchain#1{\hbox{\edgeno0%
\valign{&\de0{\the\edgeno}##\global\advance\edgeno
by 1\cr #1\cr}}}
\def\lchain#1{\vbox{\edgeno0%
\halign{&\le{\the\edgeno}0##\global\advance\edgeno 
by -1\cr #1\cr}}}

\def\turn{\ifx\chain\vchain\let\chain\hchain
\let\tbox\vbox\let\tskip\vskip\let\tss\vss
\else\let\chain\vchain\let\tbox\hbox
\let\tskip\hskip\let\tss\hss\fi}

\def\starttree{\let\chain\hchain%
\let\tbox\vbox\let\tskip\vskip\let\tss\vss%
\setbox\chainbox\vertexlabel0\smallnum0}

\def\tree#1#2{\hbox{\chain{#1}\turn\branch{#2}}}
\def\branch#1{\vbox{\halign{&\grow ##\cr #1\cr}}}
\def\grow#1#2{\ifx.#1\relax\else\tbox{%
\tskip #1\edgelength\tbox to 0 pt{\tss{\ifx\tbox\vbox%
\offinterlineskip\fi\tree #2}\tss}\tskip-#1\edgelength}\fi}

\def\sc#1{\smallcirc\uclabel{#1}}

\def\CoxeterDiagramHfour{\hbox{%
\sc1\shedge\delabel5\sc2\shedge\sc3\shedge\sc4}}